\input amstex
\documentstyle{amsppt}
\magnification=\magstep1 \NoRunningHeads

\topmatter

\title On new spectral multiplicities for ergodic maps
 \endtitle
\author  Alexandre~I.~Danilenko
\endauthor

\address
 Institute for Low Temperature Physics
\& Engineering of National Academy of Sciences of Ukraine, 47 Lenin Ave.,
 Kharkov, 61164, UKRAINE
\endaddress
\email alexandre.danilenko\@gmail.com
\endemail

\thanks
\endthanks

 \abstract It is shown that each subset of positive integers that contains 2 is realizable as the set of essential values of the multiplicity function for the Koopman operator of some weakly mixing transformation.
  \endabstract

 \endtopmatter
   \document

\head 0. Introduction
\endhead
Let $(X,\goth B,\mu)$ be a standard non-atomic probability space. Given a $\mu$-preserving (invertible) transformation $T$, we denote by $U_T$ the corresponding {\it Koopman operator} in $L^2(X,\mu)$,
$U_Tf:=f\circ T$. Let $\Cal M(T)$ stand for the set of essential values of the spectral multiplicity function for the restriction of $U_T$ to the subspace of 0-mean functions $L^2_0(X,\mu):=L^2(X,\mu)\ominus \Bbb C$.
We call a subset $M$ of positive integers {\it realizable} if there is an ergodic transformation $T$ such that $M=\Cal M(T)$. In the present paper we investigate a  long-standing open problem in the spectral theory of dynamical systems that can be stated as follows:
\roster
\item"---"
{\it What subsets of $\{1,2,\dots\}$ are realizable?}
\endroster
It is expected that all subsets are realizable. It has been already shown that all subsets containing 1 are realizable \cite{KL} (re-proved with a different argument in \cite{A2}). See also earlier works \cite{Os}, \cite{R1}, \cite{R2}, \cite{G--L} on the subject. We note that the spectral multiplicities from those papers are realized on transformations that are {\it compact group extensions of rank-one maps}.

Less is known about realizability of subsets without 1. Whether $\{n\}$ is realizable for $n>1$? This problem of Rokhlin was first solved for $n=2$ in \cite{Ag1} and \cite{Ry1}.  The  transformations considered in those papers are {\it Cartesian squares of  rank-one maps}. Other realizable sets came with n-fold Cartesian products and their natural factors: $\{n,n(n-1),\dots, n!\}$ in  \cite{Ag1}, $\{2,3,\dots,n\}$ in \cite{Ag4}, etc. It is worth to note that those works on Cartesian products were influenced by the paper \cite{Ka} which circulated since mid-eighties as an unpublished manuscript.
 As was shown in \cite{Ry2}, \cite{Ry3} and \cite{Ag4} those sets without 1 are also realizable in the class of mixing transformations.

For an arbitrary $n$, the Rokhlin problem on homogeneous spectrum was first solved  in \cite{Ag3} in a non-constructive way.  An explicit solution appeared  in \cite{Da2}. A method of {\it auxiliary non-Abelian group actions} was in use in those two papers. The explicit construction from \cite{Da2} combined with the techniques of compact group extensions was  used  to show that for each $n>1$ and a  subset $M\subset \Bbb N$, the set $n\cdot (M\cup\{1\})$ is realizable  \cite{Da2}.

 Let $G$ be a countable Abelian group, $H$ a subgroup of $G$ and $v:G\to G$ a group automorphism. We set
$$
L(G,H,v):=\{\#(\{v^i(h)\mid i\in\Bbb Z\}\cap H), h\in H\setminus\{0\}\}.
$$
It was shown in a recent paper \cite{KaL} that $\{2\}\cup L(G,H,v)$ is realizable whenever $v$ is periodic. In particular, subsets $\{2\}\cup n\cdot (M\cup\{1\})$, where $n>1$ and $M$ a finite subset of $\Bbb N$, are realizable. That answers a question from \cite{Ry3} (see also \cite{Da3, Section 5}). However, \cite{KaL} does not contribute to realization of infinite subsets because periodicity of $v$ bounds $L(G,H,v)$ to be finite. Moreover, it remains unclear whether {\it every} finite subset of $\Bbb N$ equals $L(G,H,v)$ for some triplet $(G,H,v)$ with   $v$ periodic.

The purpose of the present paper is to prove the following theorem which extends the main result of \cite{KaL} to all subsets containing 2.

\proclaim{Main Theorem} Let $E$ be an arbitrary subset of positive integers. Then there is a weakly mixing transformation $S$ such that $\Cal M(S)=E\cup\{2\}$.
\endproclaim

Our method develops further the approach considered  in \cite{KaL}. It is based upon the solution of Rokhlin problem for $n=2$ \cite{Ag1}, \cite{Ry1} and  `symmetries' of   some special compact group extensions.
We make use of the $(C,F)$-construction (see the survey \cite{Da3}) as a convenient tool to build dynamical systems and their extensions that have a prescribed `list' of {\it weak limits} for powers of Koopman operators restricted to some `components', i.e. invariant subspaces. This yields that two  components are either unitarily equivalent or spectrally disjoint. It remains to count the number of components in every unitarily equivalent class. Notice that  Katok and Lema{\' n}czyk \cite{KaL} study so-called {\it double}  (non-Abelian) compact $K\rtimes_v (\Bbb Z/n\Bbb Z)$-extensions of rank-one maps. Every such an extension can be considered as  a $K$-extension of a $\Bbb Z/n\Bbb Z$-extension, that is what `double' means here. A benefit of a double extension is that an important {\it cohomology equation} on the $K$-valued cocycle (see \thetag{2-1} below) holds automatically. In this paper we consider only {\it single} Abelian $K$-extensions. The  equation \thetag{2-1} is satisfied due to a special choice of the cocycle. An advantage of our approach is that the automorphism $v$ entering into the equation need not be periodic. This leads to realizability of infinite subsets.

\head 1. Algebraic lemma
\endhead

The following algebraic statement is a key ingredient  in the proof of  Main Theorem.

\proclaim{Algebraic Lemma} Given any subset $E\subset\Bbb N$, there exist
a countable Abelian group $G$, a subgroup $H\subset G$ and an automorphism $v:G\to G$ such that $E=L(G,H,v)$.
Moreover, the following  properties  are satisfied:
\roster
\item"\rom{(i)}"
the subgroup $\Cal G:=\{a\in\widehat G\mid a\circ v^{m_a}= a\text{ for some }m_a>0\}$ is locally finite, countable and dense in $\widehat G$,
\item"\rom{(ii)}"
if $g_1,g_2\in G$ and $v^i(g_1)\ne g_2$ for all $i\in\Bbb Z$ then there is $a\in \Cal G$  such that  and
$\sum_{i=0}^{m_a-1} a(v^i(g_1))\ne\sum_{i=0}^{m_a-1} a(v^i(g_2))$,
\item"\rom{(iii)}"
$\#\{m_a^{-1}\sum_{i=0}^{m_a-1} a(v^i(g))\mid a\in\Cal G\}=\infty$ for each $g\ne 0$.
\endroster
\endproclaim

\demo{Proof} Since the case $E=\{1\}$ is trivial, we will assume below that  $E\ne\{1\}$.
Let $n_1, n_2,\dots$ be a sequence of integers such that  $E=\{n_1,n_2,\dots,\}$ and $n_1\ne 1$. It is important that the sequence is infinite even if $E$ is finite (repetitions are allowed).
We now set  $G:=\bigoplus_{-\infty}^{+\infty}\Bbb Z/2\Bbb Z$. Let $v$ stand for the shift on $G$, i.e. if $g=(g_i)_{i\in\Bbb Z}$ then $(v(g))_i:=g_{i+1}$. To define $H$ we first construct a sequence of finite subsets $A_i\subset\Bbb Z$
such that $\# A_i\to\infty$ and
$$
2\max A_i<\min A_{i+1}
\tag 1-1
$$
for all $i=1,2,\dots$. We now let $A:=\bigsqcup_{i=1}^\infty A_i$ and set
$$
H:=\{(g_i)_{i=1}^\infty\in G\mid g_i=0\text{ if }i\notin A\text{ and }g_j=g_k\text{ whenever }j,k\in A_i\text{ for some }i\}.
$$
The subsets $A_i$ are defined via an inductive procedure. On Step $k$ we define  the subsets of cardinality $k$.

Step 1. $A_1={1}$, $A_2:=\{3\}$, \dots, $A_{n_1}:=\{3^{n_1-1}\}$.

Step $k+1$. Suppose that after Step $k$ we have already defined subsets $A_1$, \dots, $A_{l_k}$. We call a subset $B\subset\bigsqcup_{i=1}^{l_k} A_i$ \,  $(k+1)$-{\it basic} if $\# B=k+1$ and if $A_i\cap B\ne\emptyset$ for some $1\le i\le n_k$ then $A_i\subset B$. Enumerate all the $(k+1)$-basic subsets: $B_1$, \dots, $B_{r_k}$. Now we put
$$
A_{l_k+rn_{k+1}+s}:=B_{r+1}+i_{r n_{k+1}+s}, \ \text{ for all }0\le r< r_k, 1\le s< n_{k+1},
$$
where the positive integers $(i_j)_{1\le j\le (n_{k+1}-1)r_k}$ are chosen so that \thetag{1-1} is satisfied.

We now verify the conclusion of the lemma for the triplet $(G,H,v)$. Take $g=(g_i)_{i\in\Bbb Z}\in H$. Let $C:=\{i\mid g_i\ne 0\}$.
Notice that  $g_i=g_j=1$ for all $i,j\in C$. Denote by $p$ the cardinality of $C$.
Then $p$ is the smallest number such that $C\subset\bigsqcup_{i\le l_p} A_i$. By the construction, there is a $k$-basic subset $B$ such that $C$ is a translation of $B$ in $\Bbb Z$.  Moreover, there exist exactly $n_p$ different  translations of $C$ which are inside $A$. This means that the $v$-orbit of $g$ intersects $H$ exactly $n_p$ times.
 Therefore $L(G,H,v)\subset E$.
The converse inclusion is obvious. Thus the first claim of the lemma is shown.

It is easy to see that (i) is satisfied.

Let $\Cal G_m:=\{a\in\Cal G\mid \widehat v^m(a)=a\}$.
Then $\Cal G_m$ is a finite $\widehat v$-invariant subgroup of $\Cal G$.
If $g_1,g_2\in G$ and $v^i(g_1)\ne g_2$ for all $i\in\Bbb Z$ then we can find $m>0$ such that
$v^i(g_1)\restriction\Cal G_m\ne g_2\restriction\Cal G_m$ for all $i\in\Bbb Z$.
Then
$$
\bigg(\sum_{i=0}^{m} v^i(g_1)\restriction\Cal G_m\bigg)\perp\bigg(\sum_{i=0}^{m} v^i(g_2)\restriction\Cal G_m\bigg) \quad\text{as elements of }L^2(\Cal G_m)
$$
and (ii) follows.

 We consider an element $0\ne g\in G$ as an infinite  sequence of symbols $0$ and $1$ with finitely many, say $j$,  symbols $1$. Take a block $b\in(\Bbb Z/2\Bbb Z)^m$ consisting of one symbol $1$ and $m-1$ symbols $0$ for a very large $m$.
Then we set $a:=b^\infty\in \Cal G_m$.
It is easy to verify that
$m^{-1}\sum_{i=0}^{m-1} a(v^i(g))\mid a\in\Cal G\}=(m-2j)/m$.
This yields~(iii).
\qed
\enddemo

Notice that   stronger  versions  of the above result  have been established  in some  particular cases:
\roster
\item"$(\circ)$" If $1\in E$ then  $v$ in the statement of Algebraic Lemma can be chosen {\it quasi-periodic}, i.e. every $v$-orbit is finite. If, in addition, $E$ is finite than $E=L(G,H,v)$  for  finite dimensional toruses $G$ and $H$ and a periodic automorphism $v$ \cite{KL}.
\item"$(\circ)$" If $E=\{2\}\cup n\cdot (M\cup\{1\})$, where $n>1$ and $M$ a finite subset of $\Bbb N$, then $E=L(G,H,v)$  for   finite groups $G$ and $H$ and a periodic automorphism $v$  \cite{KaL}.
\endroster

In this connection we note that $v$ in our construction is not quasi-periodic  even for $E$ finite or $1\in E$.

\head 2.  Weak limits of powers, Cocycles, $(C,F)$-construction
\endhead

We will need two lemmata on spectral properties of some Cartesian products. For the proof we refer to \cite{Ag1}, \cite{Ry1} and \cite{KaL}.

\proclaim{Lemma 2.1 \text{(\cite{Ag1}, \cite{Ry1})}}
Let  $T$ be a weakly mixing transformation with simple spectrum. If the weak closure of powers of $U_T$ contains $0.5(I+U_T)$ then $T\times T$ has  homogeneous spectrum of multiplicity $2$ in the orthocomplement to the constants.
\endproclaim

We note that a theory of {\it linked approximation} suggested in \cite{KaS} plays an important role in the proof of the above lemma.

\proclaim{Lemma 2.2 \text{(\cite{KaL})}} Let $V_i$, $i=1,2$, be unitary operators with simple spectrum. Assume moreover that there are two sequences $(n_t)_{t>0}$ and $(m_t)_{t>0}$ such that
\roster
\item"\rom{(i)}"
$V_i^{n_t}\to 0.5(I+V_i^*)$ weakly, $i=1,2$.
\item"\rom{(ii)}"
$V_i^{m_t}\to 0.5(I+c_iV_i^*)$ weakly, $i=1,2$.
\endroster
If $c_1\ne c_2$ then $V_1\otimes V_2$ has also a simple spectrum.
\endproclaim

Let $T$ be an ergodic transformation of $(X,\mu)$. Denote by $\Cal R\subset X\times X$ the $T$-orbit equivalence relation. A Borel map $\alpha$ from $\Cal R$ to a compact group $K$ is called a {\it cocycle} of $\Cal R$ if
$$
\alpha(x,y)\alpha(y,z)=\alpha(x,z)\quad\text{for all }(x,y), (y,z)\in\Cal R.
$$
 Two cocycles $\alpha,\beta:\Cal R\to K$ are {\it cohomologous} if
$$
\alpha(x,y)=\phi(x)\beta(x,y)\phi(y)^{-1}\quad \text{at a.a. }(x,y)\in \Cal R
$$
for a Borel map $\phi:X\to K$. If a transformation $S$ commutes with $T$ (i.e. $S\in C(T)$) then a cocycle $\alpha\circ S:\Cal R\to K$ is well defined by $\alpha\circ S(x,y):=\alpha(Sx,Sy)$. The important cohomology equation on $\alpha$ mentioned in Section 0 can now be stated as follows
$$
\alpha\circ S\text{ \ is cohomologous to \ }  v\circ\alpha\tag 2-1
$$
for some $S\in C(T)$ and a group automorphism $v:K\to K$.

To prove Main Theorem we will use  the $(C,F)$-construction (see \cite{dJ}, \cite{Da1}--\cite{Da3}).  We now briefly outline its formalism.
Let two sequence $(C_n)_{n>0}$ and $(F_n)_{n\ge 0}$ of finite subsets in $\Bbb Z$ are given such that:
\roster
\item"---"
$F_n=\{0,1,\dots,h_n-1\}$, $h_0=1$, $\# C_n>1$,
\item"---" $F_n+C_{n+1}\subset F_{n+1}$,
\item"---" $(F_{n}+c)\cap (F_n+c')=\emptyset$ if $c\ne c'$, $c,c'\in C_{n+1}$,
\item"---" $\lim_{n\to\infty}\frac{h_n}{\#C_1\cdots\# C_n}<\infty$.
\endroster
Let $X_n:=F_n\times C_{n+1}\times C_{n+2}\times\cdots$. Endow this set with the (compact Polish) product topology. The  following map
$$
(f_n,c_{n+1},c_{n+2})\mapsto(f_n+c_{n+1},c_{n+2},\dots)
$$
is a topological embedding of $X_n$ into $X_{n+1}$. We now set $X:=\bigcup_{n\ge 0} X_n$ and endow it with the (locally compact Polish) inductive limit topology. Given $A\subset F_n$, we denote by $[A]_n$ the following cylinder: $\{x=(f,c_{n+1},\dots,)\in X_n\mid f\in A\}$. Then $\{[A]_n\mid A\subset F_n, n>0\}$ is the family of all compact open subsets in $X$. It forms a base of the topology on $X$.

Let  $\Cal R$ stand for the {\it tail} equivalence relation on $X$: two points $x,x'\in X$ are $\Cal R$-equivalent if there is $n>0$ such that $x=(f_n,c_{n+1},\dots),\ x'=(f_n',c_{n+1}',\dots)\in X_n$ and $c_m=c_m'$ for all $m>n$. There is only one probability (non-atomic) Borel measure $\mu$ on $X$ which is invariant (and ergodic) under $\Cal R$.

Now we define a transformation $T$ of $(X,\mu)$ by setting
$$
T(f_n,c_{n+1},\dots):=(1+f_n,c_{n+1},\dots )\text{ whenever }f_n<h_n-1,\ n>0.
$$
This formula defines $T$ partly on $X_n$. When $n\to\infty$, $T$ extends to the entire $X$ minus countably many points as a $\mu$-preserving transformation. Moreover, the $T$-orbit equivalence relation  coincides  with $\Cal R$ (on the subset where $T$ is defined). We call $T$ {\it the $(C,F)$-transformation} associated with $(C_{n+1},F_n)_{n\ge 0}$.

We recall a concept of $(C,F)$-cocycle (see \cite{Da2}). From now on, the group $K$ is assumed Abelian. Given a sequence of maps $\alpha_n:C_n\to K$, $n=1,2,\dots$, we first define a Borel cocycle $\alpha:\Cal R\cap (X_0\times X_0)\to K$ by setting
$$
\alpha(x,x'):=\sum_{n>0}(\alpha_n(c_n)-\alpha_n(c_n')),
$$
whenever $x=(0,c_1,c_2,\dots)\in X_0$,  $x'=(0,c_1',c_2',\dots)\in X_0$ and $(x,x')\in\Cal R$. To extend $\alpha$ to the entire $\Cal R$, we first define a map $\pi:X\to X_0$ as follows. Given $x\in X$, let $n$ be the least positive integer such that $x\in X_n$. Then $x=(f_n,c_{n+1},\dots)\in X_n$. We  set
$$
\pi(x):=(\underbrace{0,\dots,0}_{n+1\text{ times}}, c_{n+1}, c_{n+2},\dots)\in X_0.
$$
Of course, $(x,\pi(x))\in\Cal R$. Now for each pair $(x,y)\in\Cal R$, we let
$$
\alpha(x,y):=\alpha(\pi(x),\pi(y)).
$$
It is easy to verify that $\alpha$ is a well defined cocycle of $\Cal R$ with values in $K$. We call it {\it the $(C,F)$-cocycle associated with} $(\alpha_n)_{n=1}^\infty$.

The following statement follows from \cite{Da2, Section 4}.

\proclaim{Lemma 2.3} Let $\bar z=(z_n)_{n+1}^\infty$ be a sequence of positive reals. Suppose that
$$
\sum_{n>0}\# (C_n\triangle (C_n-z_n))/\#C_{n}<\infty.
 $$
For each $m>0$, we set
 $$
X^{\bar z}_m:=\{0,1,\dots, h_m-z_1-\cdots-z_m\}\times\prod_{n>m}(C_{n}\cap(C_n-z_n))\subset X_m.
$$
Then a transformation $S_{\bar z}$ of $(X,\mu)$ is well defined by setting
$$
S_{\bar z}(x):=(z_1+\cdots+z_m+f_m, z_{m+1}+c_{m+1}, z_{m+2}+c_{m+2},\dots)\tag2-2
$$
for all $x=(f_m,c_{m+1},c_{m+2},\dots)\in X^{\bar z}_m$, $m=1,2,\dots$.
Moreover, $S_{\bar z}$ commutes with $T$ and $T^{z_1+\cdots+z_m}\to S_{\bar z}$ as $m\to\infty$.

Now let $C_m^\circ:=\{c\in C_m\cap (C_m-z_m)\mid \alpha_m(c+z_m)=v(\alpha_m(c))\}$. If
$$
\sum_{n>0}(1-\# C_n^\circ/\#C_{n})<\infty\tag2-3
$$
then the cocycle $\alpha\circ S_{\bar z}$ is cohomologous to $v\circ\alpha$.
\endproclaim

\head 3. Proof of Main Theorem
\endhead

 By Algebraic Lemma, there exist a compact Polish Abelian group $K$,  a closed subgroup $H$ of $K$ and  a continuous automorphism $v$ of $K$ such that
$$
E=L(\widehat K,\widehat{K/H}, \widehat v).
 $$
We also assume that the  conditions  (i)-(iii) from the statement of the lemma are satisfied.
The subgroup  of $v$-periodic points in $K$ will be denoted by $\Cal K$.

We will construct some special $(C,F)$-transformation and its cocycle with values in $K$. Fix a partition
$$
\Bbb N=\bigsqcup_{a\in \Cal K} \Cal N_a\sqcup\bigsqcup_{a,b\in \Cal K}\Cal N_{a,b}
$$
of $\Bbb N$ into infinite subsets. Now we define a sequence $(C_n,h_n,z_n,\alpha_n)_{n=1}^\infty$ via an inductive procedure.
Suppose we have already constructed this sequence up to index $n$.
Consider now two cases.

{\bf [I]} If  $n+1\in \Cal N_a$ for some $a\in\Cal K$, we denote by $m_a$ the least positive period of $a$ under $v$.
Now we set
$$
\gather
z_{n+1}:=m_anh_n, \quad r_n:=n^3m_a,\\
C_{n+1}:=h_n\cdot\{0,1,\dots, r_n-1\},\\
h_{n+1}:= r_nh_n,
\endgather
$$
Let $\alpha_{n+1}:C_{n+1}\to K$ be any map satisfying the following conditions
\roster
\item"(A1)" $\alpha_{n+1}(c+z_{n+1})=v\circ\alpha_{n+1}(c)$ for all $c\in C_{n+1}\cap (C_{n+1}-z_{n+1})$,
\item"(A2)" for each $0\le i< m_a$ there is a subset $C_{n+1,i}\subset C_{n+1}$ such that
$$
\gather
C_{n+1,i}-h_n\subset C_{n+1}, \\  \alpha_{n+1}(c)=\alpha_{n+1}(c-h_n)+v^i(a)\text{ for all }c\in C_{n+1,i}\text{  and}\\
\bigg|\frac{\#C_{n+1,i}}{\#C_{n+1}}-\frac 1{m_a}\bigg|<\frac 2{nm_a}.
\endgather
$$
\endroster
{\bf [II]} If  $n+1 \in \Cal N_{a,b}$ for some $a,b\in\Cal K$,  we denote by $m_{a,b}$ the least  common positive period of $a$ and $b$ under $v$.
Now we set
$$
\gather
z_{n+1}:=m_{a,b}n(2h_n+1), \quad r_n:=2n^3m_{a,b},\\
D_{n+1}:=h_n\cdot\{0,1,\dots, nm_{a,b}-1\}\sqcup ((h_n+1)\cdot\{1,2,\dots,nm_{a,b}\}+h_n(nm_{a,b}-1)),\\
C_{n+1}:=D_{n+1}+z_{n+1}\cdot\{0,1,\dots,n^2-1\},\\
h_{n+1}:= r_nh_n+r_n/2.
\endgather
$$
Let $\alpha_{n+1}:C_{n+1}\to K$ be any map satisfying the following conditions
\roster
\item"(A3)" $\alpha_{n+1}(c+z_{n+1})=v\circ\alpha_{n+1}(c)$ for all $c\in C_{n+1}\cap (C_{n+1}-z_{n+1})$,
\item"(A4)" for each $e\in\{a,b\}$ and $0\le i< m_a$ there is a subset $C_{n+1,i}^e\subset C_{n+1}$ such that
$$
\gather
C_{n+1,i}^e-h_n-1\subset C_{n+1}, \\  \alpha_{n+1}(c)=\alpha_{n+1}(c-h_n-1)+v^i(e)\quad
\text{ for all }c\in C_{n+1,i}^e\text{  and}\\
\bigg|\frac{\#C_{n+1,i}^e}{\#C_{n+1}}-\frac 1{2m_{a,b}}\bigg|<\frac 2{nm_{a,b}}.
\endgather
$$
\endroster
Thus, $C_{n+1},h_{n+1}, z_{n+1},\alpha_{n+1}$ are completely defined.

We now let $F_n:=\{0,1,\dots,h_n-1\}$. Denote by $(X,\mu,T)$ the  $(C,F)$-transforma\-ti\-on associated with the sequence $(C_{n+1},F_n)_{n\ge 0}$. Let $\Cal R$ stand for the tail equivalence relation (or, equivalently, $T$-orbit equivalence relation) on $X$.
 Denote by $\alpha:\Cal R\to K$ the cocycle of $\Cal R$ associated with the sequence $(\alpha_n)_{n>0}$.
Let $\lambda_{K/H}$ stand for the Haar measure on $K/H$. We denote by $T_{\alpha,H}$ the following transformation of the space
$(X\times K/H,\lambda_{K/H})$:
$$
T_{\alpha,H}(x,k+H):=(Tx,\alpha(Tx,x)+k+H).
$$
{\it Our purpose} is to prove that $\Cal M(T\times T_{\alpha,H})=E\cup\{2\}$.

Since
$$
\sum_{n>0}
\frac{\#(C_n\triangle (C_n-z_n))}{\# C_n}=\sum_{n>0}\frac 2{n^2},
$$
it follows from Lemma 2.3 that a transformation $S_{ \bar z}$ of $(X,\mu)$ is well defined by the formula \thetag{2-2} and $S_{\bar z}\in C(T)$.

 It follows from (A1) and (A3) that \thetag{2-3} is satisfied. Hence by Lemma~2.3,  the cocycle $\alpha\circ S_{\bar z}$ is cohomologous to $v\circ\alpha$.

We need  more notation.
Given $a\in\Cal K$ and $\chi\in \widehat K$,  let
$$
l_\chi(a):= m_a^{-1}\sum_{i=0}^{m_a-1}\chi(v^i(a)).
$$
We also denote by $U_{T,\chi}$ the unitary operator on the space $L^2(X,\mu)$  given by
$$
U_{T,\chi}g(x)=\chi(\alpha(Tx,x))g(Tx),
$$

\proclaim{Lemma 3.1} Let  $a,b\in \Cal K$.
Then for each $\chi\in\widehat K$,
\roster
\item"\rom{(i)}"
$U_{T,\chi}^{h_n}\to l_\chi(a)\cdot I$
 as  $\Cal N_a-1\ni n\to\infty$ and
\item"\rom{(ii)}"
$U_{T,\chi}^{h_n}\to 0.5(l_\chi(a)\cdot I
+l_\chi(b)\cdot U_{T,\chi}^*)$ as  $\Cal N_{a,b}-1\ni n\to\infty$.
\endroster
\endproclaim

\demo{Proof} We show only (ii) since (i) is proved in a similar way but a bit simpler.
Let $n\in\Cal N(a,b)$.

 Take any subset $A\subset F_n$.  We note that $[A]_n=[A+C_{n+1}]_{n+1}$.
Therefore it follows from (A4) that
for each $x\in T[F_n]_n$,
$$
\align
U_{T,\chi}^{h_n}1_{[A]_n}(x)&=\sum_{e\in\{a,b\}}\sum_{i=0}^{m_{a,b}-1} \chi(\alpha(T^{h_n}x,x)) 1_{[A+C^e_{n+1,i}]_{n+1}}(T^{h_n}x)+ h(x)\\
&=\sum_{i=0}^{m_{a,b}-1} \chi(v^i(a)) 1_{[A+C^a_{n+1,i}-h_n]_{n+1}}(x)\\
&+\sum_{i=0}^{m_{a,b}-1} \chi(v^i(b)\alpha(T^{-1}x,x)) 1_{[A+C^b_{n+1,i}-h_n-1]_{n+1}}(T^{-1}x)+h(x)\\
&=\sum_{i=0}^{m_{a,b}-1} \chi(v^i(a)) 1_{[A+C^a_{n+1,i}-h_n]_{n+1}}(x)\\
&+\sum_{i=0}^{m_{a,b}-1} \chi(v^i(b))U_{T,\chi}^* 1_{[A+C^b_{n+1,i}-h_n-1]_{n+1}}(x)+h(x),\\
\endalign
$$
where $x\mapsto h(x)$ is a function whose $L^2$-norm is small.
Hence
$$
U_{T,\chi}^{h_n}-\sum_{i=0}^{m_{a,b}-1} \chi(v^i(a)) 1_{[C^a_{n+1,i}-h_n]_{n+1}}-\sum_{i=0}^{m_{a,b}-1} \chi(v^i(b))U_{T,\chi}^* 1_{[C^b_{n+1,i}-h_n-1]_{n+1}}\to 0
$$
weakly as $\Cal N_{a,b}-1\ni n\to\infty$, where the  functions
$$
1_{[C^a_{n+1,i}-h_n]_{n+1}}, 1_{[C^b_{n+1,i}-h_n-1]_{n+1}}\in L^\infty(X,\mu)
$$
are considered as multiplication operators in $L^2(X,\mu)$.

It remains to use the inequalities from (A2) and (A4) and a standard fact that for any sequence $C_{n}'\subset C_{n}$ such that $\# C_{n}'/\# C_n\to\delta$  for some $\delta>0$ then
$$
1_{[C_{n}']_{n}}\to \delta I \quad\text{weakly as }
n\to\infty.
$$
 \qed
\enddemo

\demo{Proof of Main Theorem}
 Lemma~3.1(ii) yields  that
$U_{T_\alpha}^{h_n}\to 0.5(I+U_{T,\alpha}^*)$ weakly as ${\Cal N}_{0,0}-1\ni n\to\infty$.
It follows that the transformation $T_\alpha$ (and hence its factor $T_{\alpha,H}$) is weakly mixing.

To show that  $\Cal M(T\times T_{\alpha,H})=E\cup\{2\}$ we consider
 a natural decomposition
of   $U_{T\times T_{\alpha,H}}$ into
 an orthogonal sum
$$
U_{T\times T_{\alpha,H}}=\bigoplus_{\chi\in \widehat {K/H}}(U_T\otimes U_{T,\chi}).
$$
It is enough
 to prove the following:
\roster
\item"(a)"
$U_T\otimes U_T$ has homogeneous spectrum 2 in the orthocomplement to the constants,
\item"(b)"
$U_T\otimes U_{T,\chi}$ has simple spectrum if $\chi\ne 0$,
\item"(c)"
$U_T\otimes U_{T,\chi}$ and $U_T\otimes U_{T,\xi}$ are unitarily equivalent if $\chi$ and $\xi$ belong to the same $\widehat v$-orbit,
\item"(d)" the measures of maximal spectral type of $U_T\otimes U_{T,\chi}$ and $U_{T}\otimes U_{T,\xi}$ are mutually singular if $\chi$ and $\xi$ are on not on the same $\widehat v$-orbit.
\endroster

By Lemma 3.1(ii), $U_T^{h_n}\to 0.5(I+U_T^*)$ as $\Cal N_{0,0}-1\ni n\to\infty$. Therefore (a) follows from Lemma~2.1.

Since $T$ is of rank one and the map $[f]\ni x\mapsto \alpha(Tx,x)\in K$ is constant for each $f\in F_n\setminus\{h_n-1\}$, $n\in\Bbb N$, it follows that the operator $U_{T,\chi}$ has simple spectrum.
Moreover,  $U_{T,\chi}^{h_n}\to 0.5(I+U_{T,\chi}^*)$ as $\Cal N_{0,0}-1\ni n\to\infty$ by Lemma~3.1(ii).
Since $\chi$ is nontrivial, it follows from  claim (ii) of  Algebraic Lemma that
 there is $a\in\Cal K$ with $l_\chi(a)\ne 1$.  By Lemma~3.1(ii),
$$
U_T^{h_n}\to 0.5(I+U_T^*)\text{ \  but \ }U_{T,\chi}^{h_n}\to 0.5(I+l_\chi(a)U_{T,\chi}^*)
 $$
as $\Cal N_{0,a}-1\ni n\to\infty$. Therefore Lemma 2.2 implies (b).

Since the cocycles $\alpha\circ S_{\bar z}$ and $v\circ\alpha$ are cohomologous, $U_{T,\chi}$ and $U_{T,\xi}$ are unitarily equivalent whenever $\chi$ and $\xi$ lie on the same orbit of $\widehat v$ (see \cite{G--L}, \cite{KL}).  This yields~(c).

To prove (d), we first find $a\in\Cal G$ such that $l_\chi(a)\ne l_\xi(a)$ (see  claim (ii) of Algebraic Lemma).
It follows from Lemma~3.1(i) that
$$
 U_T^{h_n}\otimes U_{T,\chi}^{h_n}\to l_\chi(a)I \text{ \ and \  } U_T^{h_n}\otimes U_{T,\xi}^{h_n}\to l_\xi(a)I
$$
as $\Cal N_a-1\ni n\to\infty$. Hence (d) holds. \qed
\enddemo

\head 4. Concluding remarks
\endhead

Combining the techniques developed in \cite{KaL} with our approach one can also realize some  subsets of $\Bbb N\setminus\{1,2\}$. For instance, given any subset $E\subset\Bbb N$, let $T_{\alpha,H}$ denote the skew product transformation constructed in Section 3.
Then
$$
\Cal M(T^{\times k}\times T_{\alpha,H})=\{k+1,(k+1)k,\dots, (k+1)!\}\cup \{k, k(k-1),\dots, k!\}\cdot E.\tag4-1
$$
 For example, taking $k=2$ we obtain $\{3,6\}\cup 2\cdot E$.
To show \thetag{4-1} we repeat the proof of Main Theorem almost literally but apply the following lemma instead of~Lemma~2.2.

\proclaim{Lemma 4.2 (\cite{KaL, Proposition 1}, \cite{DaR, Lemma~1.2})}
Let $V$ and $W$ be unitary operators with simple spectrum in Hilbert spaces $\Cal H$ and $\widetilde{\Cal H}$ respectively. Assume moreover that for each $i=1,\dots,k$, there are two sequences $n_t^{(i)}\to\infty$ and $m_t^{(i)}\to\infty$ and complex numbers $\kappa_i\ne \widetilde\kappa_i$ such that
\roster
\item"\rom{(i)}"
$V^{n_t^{(i)}}\to 0.5(\kappa_iI+V^*)$, $W^{n_t^{(i)}}\to 0.5(\kappa_iI+W^*)$ weakly,
\item"\rom{(ii)}"
$V^{m_t^{(i)}}\to 0.5(\kappa_iI+V_i^*)$,
$W^{m_t^{(i)}}\to 0.5(\widetilde\kappa_iI+W^*)$ weakly
\endroster
and $\#\{\kappa_1,\dots,\kappa_k\}=k$.
Then $V^{\odot k}\otimes W$ has a simple spectrum.
\endproclaim

A more general class of multiplicities arises when considering {\it natural} factors of $T^{\times k}\times T_{\alpha,H}$ as in \cite{Ag4}. Moreover, Ryzhikov recently constructed a new series of realizable sets that do not belong to this class \cite{Ry4}. For instance, sets $\{p,q,pq\}$ are realizable for all $p,q>0$.
Sets $\{3,4\}$ and $\{3,5\}$ are the simplest ones that are  not covered by this class nor  by Ryzhikov's one.

We also note that while our approach  based on Albebraic Lemma and weak limits techniques is not sufficient to solve the spectral multiplicity problem completely in the framework of ergodic {\it finite} measure preserving transformations, it can be adapted to the framework of {\it infinite} measure preserving maps to show the following: {\it every} subset of $\Bbb N$ is realizable on an ergodic infinite measure preserving transformation. This is done in a subsequent paper of Ryzhikov and the  author \cite{DaR}.

It is also interesting to study the spectral multiplicity problem for general Abelian group actions. Some results in this field were recently obtained by Konev, a post-graduate student of Ryzhikov, for $\Bbb Z^2$-actions \cite{Ko}.

\enddocument